\documentclass[11pt]{article}
\usepackage[a4paper,margin=29mm]{geometry}
\usepackage[T1]{fontenc}
\usepackage{lmodern}
\usepackage{amsmath,amssymb,amsthm,mathtools}
\usepackage{booktabs}
\usepackage{microtype}
\usepackage{tikz}
\usepackage[hidelinks]{hyperref}
\newcommand{\PaperAuthor}{Obinna Okechukwu}
\hypersetup{
  pdftitle={Graceful supertrees with at most 2m-diam(T) edges},
  pdfauthor={\PaperAuthor},
  pdfsubject={An explicit completion of a level labeling to a graceful tree},
  pdfkeywords={graceful labeling, tree embedding, graceful supertree, range-relaxed labeling}
}
\newtheorem{theorem}{Theorem}[section]
\newtheorem{lemma}[theorem]{Lemma}
\newtheorem{proposition}[theorem]{Proposition}
\newtheorem{corollary}[theorem]{Corollary}
\theoremstyle{remark}
\newtheorem{remark}[theorem]{Remark}
\newcommand{\ecc}{\operatorname{ecc}}
\newcommand{\diam}{\operatorname{diam}}

\newcommand{\interval}[2]{[#1,#2]_{\mathbb Z}}
\title{Graceful supertrees with at most\\
$2m-\diam(T)$ edges}
\author{\PaperAuthor}
\date{}

\begin{document}
\maketitle
\begin{abstract}
A graceful labeling of a tree with $m$ edges is a bijection from its vertices to $\{0,1,\ldots,m\}$ whose edge differences are pairwise distinct. The Graceful Tree Conjecture asserts that every tree admits such a labeling; equivalently, the smallest graceful tree containing a given $m$-edge tree has $m$ edges. For this embedding problem, Sethuraman, Ragukumar and Slater proved that fewer than $4m$ edges suffice, while Sethuraman and Ragukumar observed that in practice fewer than $2m$ appear to suffice. We prove the latter bound for every tree. For any root of height $h$ with $t_h$ vertices at maximum distance, we construct a graceful tree containing the original tree with exactly $M=2m-h-t_h+1$ edges; in particular, $2m-\operatorname{diam}(T)$ edges suffice. The construction completes a cyclic reflection of van Bussel's level labeling without increasing its span; each original edge difference $d$ is replaced by $M+1-d$. Unused vertex labels and missing edge differences are indexed by the same intermediate-level positions, allowing an explicit attachment rule that preserves the tree property. The best root for this construction is computable in linear time. Because every graceful supertree labeling restricts to a range-relaxed labeling, a uniform supertree bound of $(2-\varepsilon)m+O(1)$ would imply such a range-relaxed bound for all trees.
\end{abstract}



\section{Introduction}

A \emph{graceful labeling} of a tree $G$ with $M$ edges is a bijection
$F:V(G)\to\{0,1,\ldots,M\}$ for which the $M$ edge differences
$|F(u)-F(v)|$ are pairwise distinct. These differences are then precisely
$1,\ldots,M$. Introduced by Rosa as $\beta$-valuations \cite{Rosa67} and
named by Golomb \cite{Golomb72}, graceful labelings were motivated by
Ringel's tree-decomposition problem \cite{Ringel64}. Ringel's conjecture
has been proved for sufficiently large trees by methods not requiring
that the trees be graceful \cite{MPS21,KS25}. The assertion that every
tree itself is graceful remains open; see Gallian's survey \cite{Gallian}.

Instead of labeling $T$ gracefully, one can embed it as a subtree of a
graceful tree. Write $s(T)$ for the minimum number of edges of such a
supertree. The embedding problem goes back to work of Acharya
\cite{Acharya82} and Acharya, Rao and Arumugam \cite{ARA08}; related
packing and embedding results appear in \cite{SE05,SM20}. Sethuraman
and Ragukumar supplied an algorithmic construction \cite{SR14} and later
observed that in practice graceful host trees with fewer than $2m$
edges appear to suffice for an $m$-edge tree \cite{SR2015}.
A bound of the form $km$ with
$k<\lceil m/4\rceil$ was established by Sethuraman, Ragukumar and
Slater \cite{SRS2016}; they subsequently proved $s(T)<4m$
\cite{SRS2017}. We prove a uniform bound below $2m$ and a
stronger rooted form.

\begin{theorem}\label{thm:main}
Let \(T\) be a finite tree with \(m\geq 1\) edges, rooted at a vertex \(r\). Let \(h=\ecc_T(r)\), and let
\[
t_h=\bigl|\{v\in V(T):d_T(r,v)=h\}\bigr|.
\]
Then \(T\) is a subtree of a graceful tree with exactly
\begin{equation}\label{eq:main-span}
M=2m-h-t_h+1
\end{equation}
edges. Given the rooted tree and an order on the children of each vertex, the tree and its graceful labeling can be constructed using \(O(m)\) arithmetic operations on integers of \(O(\log(m+1))\) bits.
\end{theorem}

If $r$ is an endpoint of a diameter of $T$, then $h=\diam(T)$ and
$t_h\ge1$. Hence Theorem~\ref{thm:main} immediately yields
$s(T)\le2m-\diam(T)$. We optimize the root in
Corollary~\ref{cor:bounds} and show that this optimization itself takes
$O(m)$ time (Lemma~\ref{lem:reroot}).

The numerical span in \eqref{eq:main-span} occurs already in the
level-labeling method of van Bussel \cite{vanBussel}. Our labeling
\eqref{eq:old-labels} is obtained from his provisional coordinates by a
cyclic reflection, not by an ordinary reversal of the integer interval;
Proposition~\ref{prop:cyclic} proves the exact relation. Sethuraman,
Ragukumar and Slater improved van Bussel's range-relaxed labeling
bound. They also proved that every $m$-edge tree is contained in a
graceful tree with fewer than $4m$ edges and is a spanning subtree of a
graceful planar graph \cite{SRS2017}. These results do not give a
graceful tree host of the shorter span in \eqref{eq:main-span}.
Sethuraman and Ragukumar also complete their own weighted level
labeling without increasing its span \cite{SRArs2017}. Here the labeling
\eqref{eq:old-labels} is completed to a graceful \emph{tree} at its
exact span.

The construction is controlled by the consecutive blocks of children in
breadth-first order. The omitted vertex labels and the omitted edge
differences turn out to have the same indices: every positive rank on an
intermediate level. Assigning each missing difference its corresponding
unused label reduces completion to a formula for a single attachment.
The resulting auxiliary edges form directed chains terminating at the
original tree. The nonnesting of such breadth-first blocks is related to
the queue-layout viewpoint of Heath and Rosenberg \cite{HR92}; the exact
arithmetic partition and the attachments are proved below. The labeling
is bipartite in the sense that all even-level labels lie below all
odd-level labels, a familiar feature of bipartite labeling methods
\cite{RS95}.

The next section describes the labeling and its two complementary
partitions. Section~\ref{sec:completion} constructs the host tree, and
Section~\ref{sec:consequences} derives bounds and discusses the
range-relaxed problem.

\section{The level labeling and its missing values}\label{sec:labeling}

Throughout, \(\interval{a}{b}=\{a,a+1,\ldots,b\}\), with the convention that this set is empty when \(a>b\). Fix a root \(r\) of \(T\), and order the vertices at distance \(i\) from it as
\[
V_i=(v_{i,0},\ldots,v_{i,t_i-1}),\qquad 0\leq i\leq h.
\]
Here \(t_0=1\). Each child level is ordered in consecutive blocks following the order of its parents; within each block the order is arbitrary. Put
\begin{equation}\label{eq:offsets}
L_i=\sum_{k<i}t_k,\qquad
\pi(v_{i,j})=L_i+j.
\end{equation}
Thus \(L_0=0\), \(L_1=1\), and \(L_{h+1}=m+1\).

For \(i<h\), let \(c_{i,j}\) be the number of children of \(v_{i,j}\), and put
\[
P_{i,j}=\sum_{k<j}c_{i,k}.
\]
The children of \(v_{i,j}\) have ranks
\(\interval{P_{i,j}}{P_{i,j}+c_{i,j}-1}\) in \(V_{i+1}\). Vertices with no children are allowed at every level.

Fix \(M\) as in \eqref{eq:main-span}. Alternating the direction of the labels separates the two parity classes, while subtracting one position at every intermediate level saves \(h-1\) positions compared with the uncompressed breadth-first coordinates. Define
\begin{equation}\label{eq:old-labels}
f(v_{i,j})=
\begin{cases}
L_i+j-i/2,& i\text{ even},\\[2pt]
M+1-L_i-j+(i-1)/2,& i\text{ odd}.
\end{cases}
\end{equation}
The auxiliary vertices will be indexed by
\[
\mathcal I=\{(i,j):1\leq i<h,\ 1\leq j<t_i\}.
\]
For these indices define
\begin{equation}\label{eq:new-labels}
g_{i,j}=
\begin{cases}
M+1-L_i-j+i/2,& i\text{ even},\\[2pt]
L_i+j-(i+1)/2,& i\text{ odd}.
\end{cases}
\end{equation}

\begin{lemma}\label{lem:labels}
The map \(f\) is injective, and its values together with the values \(g_{i,j}\), for \((i,j)\in\mathcal I\), partition \(\interval{0}{M}\). In particular,
\begin{equation}\label{eq:aux-count}
|\mathcal I|=\sum_{i=1}^{h-1}(t_i-1)=m-h-t_h+1=M-m.
\end{equation}
\end{lemma}

\begin{proof}
We compute the intervals occupied by even and odd levels separately.
Their internal gaps are precisely the proposed auxiliary labels; a final
endpoint calculation shows that the two families meet.

On an even level, the labels form the interval
\[
\interval{L_i-i/2}{L_i+t_i-1-i/2}.
\]
These intervals increase with depth. Between the intervals for the even levels \(i-1\) and \(i+1\), where \(i\) is odd, the unused labels are
\[
L_i+j-(i+1)/2,\qquad 1\leq j<t_i.
\]
They are exactly the corresponding values \(g_{i,j}\). On odd levels the label intervals decrease with depth. Between levels \(i-1\) and \(i+1\), where \(i\) is even, their gap consists of
\[
M+1-L_i-j+i/2,\qquad 1\leq j<t_i,
\]
again the values \(g_{i,j}\).

It remains to check that the two parity classes meet without a gap or an overlap. If \(h=2k\), the largest even-level label is \(m-k\), and the smallest odd-level label is \(m-k+1\). If \(h=2k+1\), these two labels are \(m-t_h-k\) and \(m-t_h-k+1\), respectively. These identities follow by substituting \(L_h=m+1-t_h\) and \eqref{eq:main-span} in \eqref{eq:old-labels}. The smallest label is \(f(r)=0\), and the largest is \(f(v_{1,0})=M\). Thus all the stated intervals and gaps form a disjoint partition of \(\interval{0}{M}\). Summing the sizes of the intermediate gaps gives \eqref{eq:aux-count}.
\end{proof}

The formulas in \eqref{eq:old-labels} can be compared directly with
van Bussel's provisional labeling, in which the two parity classes lie
on opposite sides of zero. The following identity describes the
conversion and the effect on every edge difference.

\begin{proposition}\label{prop:cyclic}
Define a labeling $\ell$ of $T$ by
\begin{equation}\label{eq:vanbussel-provisional}
\ell(v_{i,j})=\begin{cases}
-L_i-j+i/2,&i\text{ even},\\[2pt]
L_i+j-(i-1)/2,&i\text{ odd}.
\end{cases}
\end{equation}
These are van Bussel's provisional level formulas in breadth-first
coordinates, extended to any order with consecutive child blocks;
van Bussel places a longest path from the root leftmost. The span of
$\ell$ is $M$. Moreover,
\[
 f(v)\equiv-\ell(v)\pmod{M+1},\qquad 0\le f(v)\le M,
\]
and, for every edge $uv$ of $T$,
\begin{equation}\label{eq:complement-differences}
 |f(u)-f(v)|=M+1-|\ell(u)-\ell(v)|.
\end{equation}
\end{proposition}

\begin{proof}
Van Bussel assigns the first vertex of level $i$ provisional label
$-L_i+i/2$ for even $i$ and $L_i-(i-1)/2$ for odd $i$. Subsequent labels
decrease by one at even levels and increase by one at odd levels,
yielding \eqref{eq:vanbussel-provisional}. These formulas depend on
positions and level sizes, so they also define $\ell$ for any
consecutive-child-block order.

For even $i$, the formulas give $\ell(v_{i,j})=-f(v_{i,j})\le0$.
For odd $i$, they give
$\ell(v_{i,j})=M+1-f(v_{i,j})\ge1$.
Thus $f$ is the least nonnegative residue of $-\ell$ modulo $M+1$.
Let $a$ be the largest even-level label of $f$. The proof of
Lemma~\ref{lem:labels} shows that the smallest odd-level label is
$a+1$. Hence $\min\ell=-a$ and
$\max\ell=M+1-(a+1)=M-a$, so $\ell$ has span $M$.
Each edge has an even-level endpoint $u$ and an odd-level endpoint
$v$. Since $\ell(u)\le0<\ell(v)$ and
$f(u)<f(v)$, we have
\[
 f(v)-f(u)=M+1-\ell(v)+\ell(u)
 =M+1-|\ell(u)-\ell(v)|,
\]
which is \eqref{eq:complement-differences}.
\end{proof}

\begin{remark}\label{rem:fixed-labeling-optimal}
The labels $0$ and $M$ already occur in $f$. Thus any graceful tree
containing this \emph{labeled} copy of $T$ must have at least $M$ edges.
The construction in Theorem~\ref{thm:main} attains this lower bound for
extensions of this particular labeling. This does not imply that $M$ is
the minimum possible supertree size $s(T)$.
\end{remark}

The same intermediate ranks index the missing edge differences. To see this, it is convenient to subtract the transition depth from each endpoint sum. Define
\begin{equation}\label{eq:starts}
A_i=L_i+L_{i+1}-i\quad(0\leq i<h),
\qquad A_h=M+1,
\end{equation}
and, for \((i,j)\in\mathcal I\), put
\begin{equation}\label{eq:missing}
q_{i,j}=A_i+P_{i,j}+j-1,
\qquad d_{i,j}=M+1-q_{i,j}.
\end{equation}

\begin{lemma}\label{lem:differences}
The edge differences of \(f\) are positive and pairwise distinct. Their complement in \(\interval{1}{M}\) consists exactly of the pairwise distinct integers \(d_{i,j}\), for \((i,j)\in\mathcal I\).
\end{lemma}

\begin{proof}
For each level transition, the adjusted endpoint sums of a parent's
children form a consecutive interval. Successive block starts are
separated by the block length plus one, giving exactly one omitted
integer per positive parent rank.

For an edge \(uv\) between levels \(i\) and \(i+1\), the even-level endpoint has the smaller label by Lemma~\ref{lem:labels}. Formula~\eqref{eq:old-labels} gives
\begin{equation}\label{eq:edge-difference}
|f(u)-f(v)|=M+1-\bigl(\pi(u)+\pi(v)-i\bigr).
\end{equation}
The adjusted sums for the children of \(v_{i,j}\) form the interval
\begin{equation}\label{eq:block}
\interval{A_i+P_{i,j}+j}{A_i+P_{i,j}+j+c_{i,j}-1}.
\end{equation}
The next interval starts two integers after the upper endpoint of \eqref{eq:block}, so the intervening integer is \(A_i+P_{i,j+1}+j\). Empty child blocks cause no exception: they contribute no edge sum and still leave their one intervening integer. More explicitly, if
\(b_j=A_i+P_{i,j}+j\), then
\[
b_{j+1}=b_j+c_{i,j}+1,
\]
so the intervals \(\interval{b_j}{b_j+c_{i,j}-1}\), together with
\(b_j-1\) for \(1\leq j<t_i\), partition
\[
\interval{A_i}{A_i+t_i+t_{i+1}-2}.
\]
For every \(0\leq i<h\), the definitions give
\[
A_{i+1}-A_i=t_i+t_{i+1}-1.
\]
For \(i=h-1\), this uses \(L_h=m+1-t_h\) and the value of \(M\).
Hence the last displayed interval is \(\interval{A_i}{A_{i+1}-1}\). Since \(A_0=1\) and \(A_h=M+1\), the transition intervals partition \(\interval{1}{M}\). At level zero there is no positive rank. The omitted adjusted sums are therefore exactly \(q_{i,j}\) for \((i,j)\in\mathcal I\). Reflecting this partition by \(q\mapsto M+1-q\) proves the assertion.
\end{proof}

\section{Completion inside the same interval}\label{sec:completion}

We now add one vertex \(x_{i,j}\) carrying label \(g_{i,j}\) for each \((i,j)\in\mathcal I\). Its first index records the level of \(v_{i,j}\) in the original tree. The desired difference \(d_{i,j}\) determines a useful host label:
\begin{equation}\label{eq:host-label}
y_i(P)=
\begin{cases}
L_{i+1}+P-i/2-1,&i\text{ even},\\[2pt]
M+1-L_{i+1}-P+(i+1)/2,&i\text{ odd}.
\end{cases}
\end{equation}
The table below locates a vertex carrying that label. In its second row the host is another auxiliary vertex whose first index is one larger.

\begin{lemma}\label{lem:attachments}
For every \((i,j)\in\mathcal I\), put \(P=P_{i,j}\) and choose the host \(w_{i,j}\) as follows:
\begin{center}
\begin{tabular}{@{}ll@{}}
\toprule
Condition & Host \(w_{i,j}\)\\
\midrule
\(P=0\) & \(v_{i,t_i-1}\)\\[2pt]
\(i<h-1,\quad 0<P<t_{i+1}\) & \(x_{i+1,P}\)\\[2pt]
\(i<h-1,\quad P=t_{i+1}\) & \(v_{i+2,0}\)\\[2pt]
\(i=h-1,\quad P>0\) & \(v_{h,t_h-P}\)\\
\bottomrule
\end{tabular}
\end{center}
Each host exists and has label \(y_i(P)\). The edge \(x_{i,j}w_{i,j}\) has difference \(d_{i,j}\).
\end{lemma}

\begin{proof}
We select the host by its required label so that all four cases reduce to one subtraction. The cases cover every possibility because \(0\leq P\leq t_{i+1}\).

If \(P=0\), substitution in \eqref{eq:old-labels} gives
\(f(v_{i,t_i-1})=y_i(0)\).
If \(i<h-1\) and \(0<P<t_{i+1}\), the index \((i+1,P)\) belongs to \(\mathcal I\), and \eqref{eq:new-labels} gives
\(g_{i+1,P}=y_i(P)\).
If \(i<h-1\) and \(P=t_{i+1}\), then \(L_{i+2}=L_{i+1}+P\), so \eqref{eq:old-labels} gives
\(f(v_{i+2,0})=y_i(P)\).

For the remaining case, \(i=h-1\) and \(1\leq P\leq t_h\), so the rank \(t_h-P\) belongs to \(\interval{0}{t_h-1}\). If \(i=2k\), then
\[
M=2L_h+t_h-2k-2,
\]
and
\[
f(v_{h,t_h-P})
=M+1-L_h-t_h+P+k
=L_h+P-k-1
=y_i(P).
\]
If \(i=2k+1\), then
\[
M=2L_h+t_h-2k-3,
\]
and
\[
f(v_{h,t_h-P})
=L_h+t_h-P-k-1
=M+1-L_h-P+k+1
=y_i(P).
\]
Thus the last row also supplies the prescribed label.

Finally, \eqref{eq:new-labels} and \eqref{eq:host-label} give
\[
\begin{cases}
g_{i,j}-y_i(P)=M+2-L_i-L_{i+1}-j-P+i,&i\text{ even},\\
y_i(P)-g_{i,j}=M+2-L_i-L_{i+1}-j-P+i,&i\text{ odd}.
\end{cases}
\]
The right side is \(d_{i,j}\) by \eqref{eq:starts} and \eqref{eq:missing}, and it is positive by Lemma~\ref{lem:differences}. Hence the added edge has the required difference.
\end{proof}

\begin{proof}[Proof of Theorem~\ref{thm:main}]
The two complement identities supply one auxiliary vertex for each missing difference. We orient each prescribed edge toward an original vertex or an auxiliary vertex with larger first index, so that filling these values preserves the tree structure.

Start with \(T\), add all vertices \(x_{i,j}\), and add the edges \(x_{i,j}w_{i,j}\) of Lemma~\ref{lem:attachments}. Each added edge has a new endpoint, so every original edge is retained and no edge is added between two original vertices. No added edge is a loop. Two added edges cannot coincide: each has a different auxiliary tail, and an edge between auxiliary vertices always joins indices \(i\) and \(i+1\), with its tail at the smaller index.

Orient every added edge from \(x_{i,j}\) toward its host. If the host is auxiliary, its first index is \(i+1\). Repeating this step must reach an original vertex, since the largest auxiliary index is \(h-1\). Therefore the resulting graph \(\widehat T\) is connected. By \eqref{eq:aux-count}, it has \(M+1\) vertices and \(M\) edges, so it is a tree.

The original and auxiliary labels are exactly \(\interval{0}{M}\) by Lemma~\ref{lem:labels}. Lemmas~\ref{lem:differences} and~\ref{lem:attachments} show that the original and added edge differences are disjoint and together give \(\interval{1}{M}\), each once. This labels \(\widehat T\) gracefully.

For the complexity bound, breadth-first traversal computes the levels and their orders. A scan of the child lists computes all \(c_{i,j}\), \(P_{i,j}\), and \(L_i\). The formulas then create one label and at most one added edge per intermediate vertex. Every arithmetic operation uses integers bounded by a constant multiple of \(m\).
\end{proof}

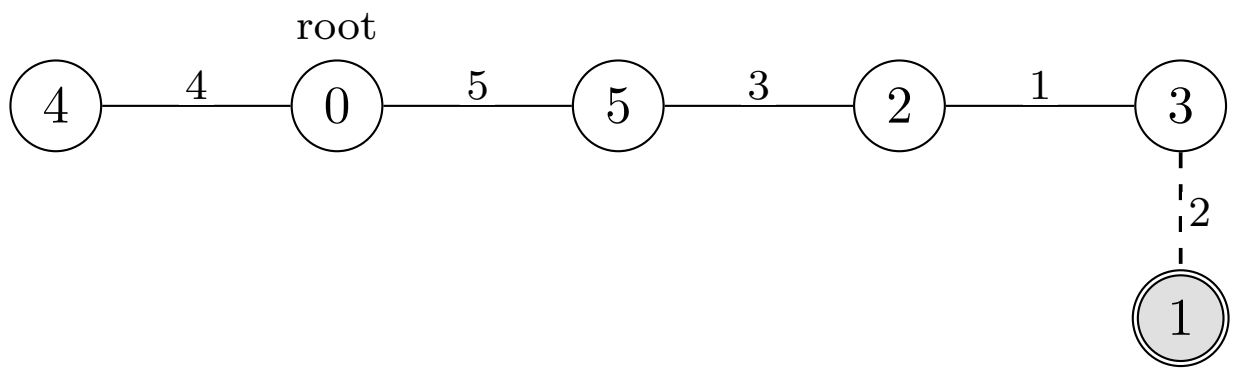
\begin{figure}[tbp]
\centering
\begin{tikzpicture}[
 x=1.24cm,y=1.04cm,
 original vertex/.style={circle,draw,fill=white,minimum size=6mm,inner sep=1pt,font=\small},
 added vertex/.style={original vertex,fill=black!12,double,double distance=.55pt},
 difference/.style={font=\scriptsize,fill=white,inner sep=1.1pt}]
 \node[original vertex] (p4) at (0,0) {4};
 \node[original vertex,label={above=6pt:{\scriptsize root}}] (p0) at (1.5,0) {0};
 \node[original vertex] (p5) at (3,0) {5};
 \node[original vertex] (p2) at (4.5,0) {2};
 \node[original vertex] (p3) at (6,0) {3};
 \node[added vertex] (p1) at (6,-1.35) {1};
 \draw (p4)--node[difference,above] {$4$} (p0);
 \draw (p0)--node[difference,above] {$5$} (p5);
 \draw (p5)--node[difference,above] {$3$} (p2);
 \draw (p2)--node[difference,above] {$1$} (p3);
 \draw[dashed,line width=.7pt] (p3)--node[difference,right] {$2$} (p1);
\end{tikzpicture}
\caption{The five-vertex path with compressed labels $4,0,5,2,3$ and
its graceful completion. The dashed edge joins the unused label $1$ to
$3$ and supplies the only missing difference, $2$.}
\label{fig:short-example}
\end{figure}

\begin{remark}[Comparison with the earlier greedy completion]\label{rem:greedy}
The in-span attachment rule of Sethuraman and Ragukumar
\cite[Graceful Tree Embedding Algorithm, Steps~2--4]{SRArs2017}
uses the smallest missing vertex label $a$. Let $V$ denote the set
of original vertex labels and $E$ the set of original edge differences.
In the algorithm's second case, $a\in E\setminus V$ and
$b=\min((V\setminus\{0\})\setminus E)$; a new vertex labeled $a$
is attached to the vertex labeled $a-b$. This one-sided rule cannot
be transferred unchanged to \eqref{eq:old-labels}.
For the five-vertex path rooted at its second vertex, our labeling is
$4,0,5,2,3$ along the path (Figure~\ref{fig:short-example}). Its edge
differences are $4,5,3,1$, so $a=1$ and $b=2$. The indicated host
$a-b=-1$ is absent. Lemma~\ref{lem:attachments} instead attaches the
new vertex $1$ to the old vertex $3$, supplying difference $2$.
The correctness proof of the algorithm in \cite{SRArs2017} uses its
\emph{own} input labeling; the example only shows why its rule cannot
be applied directly to the shorter labeling used here.
\end{remark}

\begin{figure}[tbp]
\centering
\begin{minipage}[t]{.48\textwidth}
\centering
\begin{tikzpicture}[
  x=1.2cm,y=.9cm,
  original vertex/.style={circle,draw,fill=white,minimum size=5.1mm,inner sep=.7pt,font=\small},
  difference/.style={font=\scriptsize,fill=white,inner sep=.7pt},
  every path/.style={line width=.45pt}]
  \node[original vertex] (a0) at (0,0) {0};
  \node[original vertex] (a1) at (-2.1,-1.15) {11};
  \node[original vertex] (a2) at (-.7,-1.15) {10};
  \node[original vertex] (a3) at (.7,-1.15) {9};
  \node[original vertex] (a4) at (2.1,-1.15) {8};
  \node[original vertex] (a5) at (-.7,-2.4) {4};
  \node[original vertex] (a6) at (.7,-2.4) {5};
  \node[original vertex] (a7) at (-.7,-3.65) {6};
  \draw (a0)--node[difference,above left,pos=.56] {$11$} (a1);
  \draw (a0)--node[difference,left,pos=.55] {$10$} (a2);
  \draw (a0)--node[difference,right,pos=.55] {$9$} (a3);
  \draw (a0)--node[difference,above right,pos=.56] {$8$} (a4);
  \draw (a2)--node[difference,left] {$6$} (a5);
  \draw (a3)--node[difference,right] {$4$} (a6);
  \draw (a5)--node[difference,left] {$2$} (a7);
\end{tikzpicture}
\par\smallskip
{\small (a) The labeled original tree.}
\end{minipage}\hfill
\begin{minipage}[t]{.48\textwidth}
\centering
\begin{tikzpicture}[
  x=1.2cm,y=.9cm,
  original vertex/.style={circle,draw,fill=white,minimum size=5.1mm,inner sep=.7pt,font=\small},
  added vertex/.style={original vertex,fill=black!10,double,double distance=.45pt},
  difference/.style={font=\scriptsize,fill=white,inner sep=.7pt},
  added edge/.style={dashed,line width=.6pt},
  every path/.style={line width=.45pt}]
  \node[original vertex] (b0) at (0,0) {0};
  \node[original vertex] (b1) at (-2.1,-1.15) {11};
  \node[original vertex] (b2) at (-.7,-1.15) {10};
  \node[original vertex] (b3) at (.7,-1.15) {9};
  \node[original vertex] (b4) at (2.1,-1.15) {8};
  \node[original vertex] (b5) at (-.7,-2.4) {4};
  \node[original vertex] (b6) at (.7,-2.4) {5};
  \node[original vertex] (b7) at (-.7,-3.65) {6};
  \node[added vertex] (x1) at (2.1,-2.4) {1};
  \node[added vertex] (x2) at (2.1,-3.65) {2};
  \node[added vertex] (x3) at (-2.1,-3.65) {3};
  \node[added vertex] (x7) at (.7,-3.65) {7};
  \draw (b0)--node[difference,above left,pos=.56] {$11$} (b1);
  \draw (b0)--node[difference,left,pos=.55] {$10$} (b2);
  \draw (b0)--node[difference,right,pos=.55] {$9$} (b3);
  \draw (b0)--node[difference,above right,pos=.56] {$8$} (b4);
  \draw (b2)--node[difference,left] {$6$} (b5);
  \draw (b3)--node[difference,right] {$4$} (b6);
  \draw (b5)--node[difference,left] {$2$} (b7);
  \draw[added edge] (b4)--node[difference,right] {$7$} (x1);
  \draw[added edge] (x2)--node[difference,above] {$5$} (x7);
  \draw[added edge] (x3)--node[difference,above] {$3$} (b7);
  \draw[added edge] (x7)--node[difference,above] {$1$} (b7);
\end{tikzpicture}
\par\smallskip
{\small (b) The graceful completion.}
\end{minipage}
\caption{A rooted tree with $m=7$, height $h=3$, and one vertex on the last level. The completed tree has $M=11$ edges. Vertex labels are inside circles; edge differences are written beside the edges. Added vertices have double circles and shaded interiors, and added edges are dashed. The four added edges supply the missing differences $1,3,5,7$.}
\label{fig:completion-example}
\end{figure}
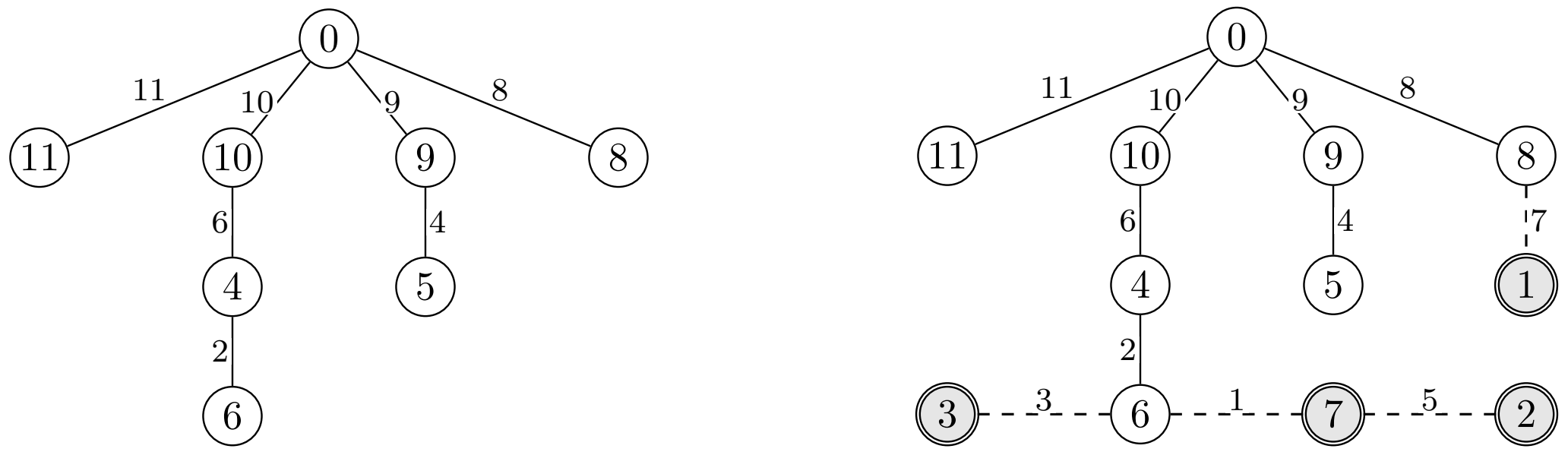

The example in Figure~\ref{fig:completion-example} uses level sizes \(1,4,2,1\). The four auxiliary vertices \(x_{1,1}\), \(x_{1,2}\), \(x_{1,3}\), and \(x_{2,1}\) have labels \(1,2,3,7\), respectively. Their prefixes are \(0,1,2,1\), so they realize all four rows of Lemma~\ref{lem:attachments}. The added edges, written by endpoint labels, are
\[
1\!-\!8,\qquad 2\!-\!7,\qquad 3\!-\!6,\qquad 7\!-\!6,
\]
with differences \(7,5,3,1\). Together with the original differences \(2,4,6,8,9,10,11\), they give \(\interval{1}{11}\).

\section{Choice of root and consequences}\label{sec:consequences}

The saving in \eqref{eq:main-span} depends on a root only through its
eccentricity and the size of its farthest level. For a vertex $r$ write
\[
 t(r)=|\{v\in V(T):d_T(r,v)=\ecc_T(r)\}|,
 \quad \zeta(T)=\max_{r\in V(T)}\bigl(\ecc_T(r)+t(r)\bigr),
\]
and put $\tau(T)=\max_r t(r)$. The parameter $\zeta(T)$ records the
largest saving obtainable by our rooted construction.

\begin{corollary}\label{cor:bounds}
For every tree $T$ with $m\ge1$ edges and diameter $D$, put
$\gamma_D(T)=\max\{t(r):\ecc_T(r)=D\}$. Then
\begin{equation}\label{eq:optimized}
 s(T)\le2m+1-\zeta(T)
 \le2m-D-\gamma_D(T)+1
 \le2m-D.
\end{equation}
In particular, $s(T)\le2m-2$ if $m\ge2$. If $T$ is not a star, then
\[
 s(T)\le2m-\tau(T)-1,\qquad s(T)\le2m-3.
\]
For a path or a star, $s(T)=m$.
\end{corollary}

\begin{proof}
Choose a root maximizing $\ecc_T(r)+t(r)$ in
Theorem~\ref{thm:main}. This proves the first inequality of
\eqref{eq:optimized}. A vertex with eccentricity $D$ lies at
an endpoint of some diameter. Choosing one with $t(r)=\gamma_D(T)$
gives $\zeta(T)\ge D+\gamma_D(T)$. Since
$\gamma_D(T)\ge1$, the other two inequalities in
\eqref{eq:optimized} follow. If $m\ge2$, then $D\ge2$.

For a nonstar, every root has eccentricity at least two. In particular,
choosing a root with $t(r)=\tau(T)$ gives $\zeta(T)\ge\tau(T)+2$,
which proves the stated $\tau$ bound. Also $D\ge3$, so the diameter
bound gives $s(T)\le2m-3$. A path rooted at an endpoint has $h=m$ and
$t_h=1$; a star rooted at its center has $h=1$ and $t_h=m$. In both
cases Theorem~\ref{thm:main} gives a host with $m$ edges. Every host
containing $T$ has at least $m$ edges.
\end{proof}

We can choose a root attaining $\zeta(T)$ without performing a search
from every vertex. The count of farthest vertices is needed as well as
their distance, so standard eccentricity rerooting must retain the
multiplicity of each maximum.

\begin{lemma}\label{lem:reroot}
For an $m$-edge tree, the values $\ecc_T(r)$ and $t(r)$ for every vertex
$r$, and hence $\zeta(T)$ and a maximizing root, can be computed in
$O(m)$ arithmetic operations.
\end{lemma}

\begin{proof}
We aggregate pairs $(d,c)$ of nonnegative integers, interpreted as a
largest distance $d$ and the number $c$ of vertices at that distance.
The operation $\oplus$ keeps the pair with larger first coordinate and
sums second coordinates in case of equality. It is associative and
commutative; an empty aggregate has a formal neutral element
$(-\infty,0)$.

For each oriented edge $u\to v$, let $B_{u\to v}$ be the pair giving
the largest distance from $v$ to a vertex in the component of $T-uv$
containing $u$, together with the number of vertices attaining it. Then
\begin{equation}\label{eq:reroot-message}
 B_{u\to v}=
 \bigl(1+d,c\bigr),\qquad
 (d,c)=(0,1)\oplus\!\!\bigoplus_{w\in N(u)\setminus\{v\}}
 B_{w\to u}.
\end{equation}
The base pair counts $u$ itself. Root the tree arbitrarily. A
postorder traversal computes every message from a child to its parent
using \eqref{eq:reroot-message}. At each vertex, prefix and suffix
aggregates of its incoming messages then compute the aggregate excluding
any specified neighbor in constant time per neighbor. A preorder
traversal consequently computes every remaining message. There are
$2m$ oriented edges, and all prefix and suffix calculations together
perform $O(\sum_v\deg(v))=O(m)$ operations.

For each vertex $r$, the aggregate
\[
 (0,1)\oplus\bigoplus_{u\in N(r)}B_{u\to r}
\]
is exactly $(\ecc_T(r),t(r))$: the components at its incident edges
partition $V(T)\setminus\{r\}$, and their farthest-distance counts
combine without overlap. Maximizing the sum of the two coordinates
finds $\zeta(T)$ and an attaining root in the same time bound.
\end{proof}

\begin{figure}[tbp]
\centering
\begin{tikzpicture}[
  x=1.13cm,y=1cm,
  original vertex/.style={circle,draw,fill=white,minimum size=5.1mm,inner sep=.7pt,font=\small},
  added vertex/.style={original vertex,fill=black!10,double,double distance=.45pt},
  difference/.style={font=\scriptsize,fill=white,inner sep=.7pt},
  added edge/.style={dashed,line width=.6pt},
  every path/.style={line width=.45pt}]
  \node[original vertex,label={above:$r$}] (c0) at (0,0) {0};
  \node[original vertex] (c9) at (1.6,0) {9};
  \node[original vertex] (c1) at (3.2,0) {1};
  \node[original vertex] (c8) at (4.8,0) {8};
  \node[original vertex] (c2) at (6.4,1.3) {2};
  \node[original vertex] (c3) at (6.4,0) {3};
  \node[original vertex] (c4) at (6.4,-1.3) {4};
  \node[original vertex] (c5) at (8,0) {5};
  \node[added vertex] (x6) at (9.6,0) {6};
  \node[added vertex] (x7) at (8,-1.3) {7};
  \draw (c0)--node[difference,above] {$9$} (c9);
  \draw (c9)--node[difference,above] {$8$} (c1);
  \draw (c1)--node[difference,above] {$7$} (c8);
  \draw (c8)--node[difference,above left] {$6$} (c2);
  \draw (c8)--node[difference,above] {$5$} (c3);
  \draw (c8)--node[difference,below left] {$4$} (c4);
  \draw (c3)--node[difference,above] {$2$} (c5);
  \draw[added edge] (c5)--node[difference,above] {$1$} (x6);
  \draw[added edge] (c4)--node[difference,above] {$3$} (x7);
\end{tikzpicture}
\caption{The tree from Figure~\ref{fig:completion-example}, rooted at an endpoint of a diameter and relabeled by the construction. Here \(m=7\), \(h=\diam(T)=5\), and the last level has one vertex. The two added vertices, labeled \(6\) and \(7\), supply the missing differences \(1\) and \(3\), giving \(M=9=2m-\diam(T)\) edges.}
\label{fig:diameter-example}
\end{figure}
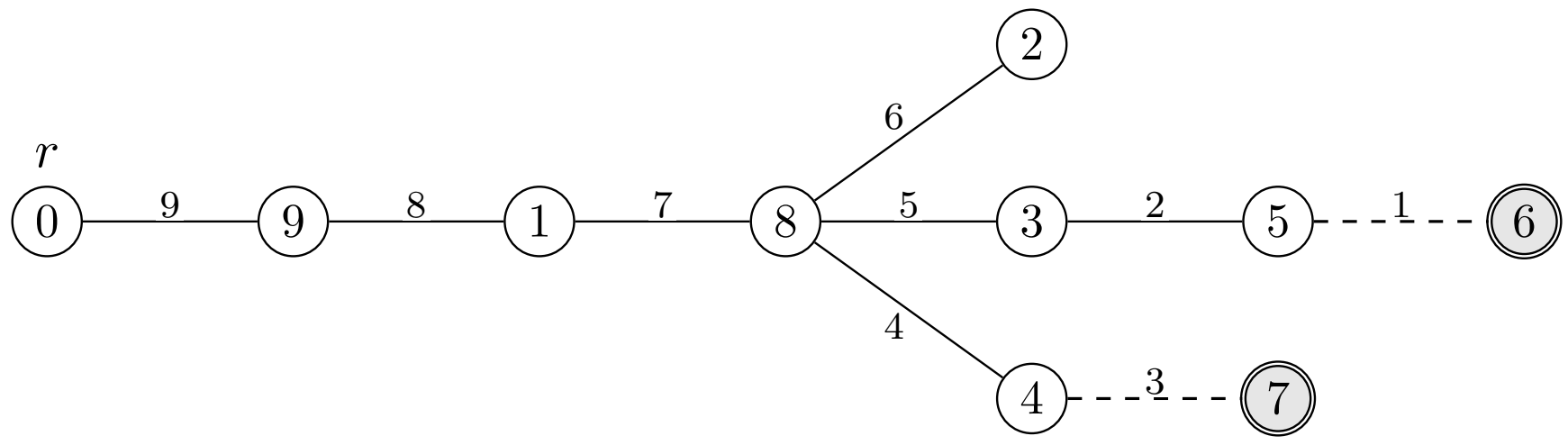

Figure~\ref{fig:diameter-example} shows the same original tree as
Figure~\ref{fig:completion-example}, now rooted at a diameter endpoint.

\section{Range-relaxed labelings and open questions}

The \emph{gracefulness} $\operatorname{grac}(T)$, as defined by
Pereira, Singh and Arumugam \cite{PSA15}, is the least integer $R$ for
which $T$ admits an injective labeling into $\{0,\ldots,R\}$ with
pairwise distinct edge differences. Such a labeling is called
\emph{range-relaxed} by van Bussel \cite{vanBussel}. Restricting a
graceful labeling of any host tree containing $T$ gives one, so
\begin{equation}\label{eq:grac-order}
 m\le\operatorname{grac}(T)\le s(T).
\end{equation}
Our construction completes the specific labeling
\eqref{eq:old-labels} without increasing its span. Equality in
\eqref{eq:grac-order} is not asserted, nor is an in-span completion
proved for arbitrary range-relaxed labelings.

Van Bussel asked whether an absolute $\varepsilon>0$ permits
range-relaxed labels of size at most $(2-\varepsilon)m$ for every
$m$-edge tree \cite{vanBussel}. Barrientos and Krop proved that
every $m$-edge tree admits range-relaxed labels of span at most
$m+\operatorname{ex}(T)$, where $\operatorname{ex}(T)$ denotes the excess
of a suitable layered representation \cite[Theorem~2.2]{BK17}. In
particular, their lobster bound is $(3m-D)/2$ for a lobster with
$m$ edges and diameter $D$ \cite[Theorem~2.4]{BK17}.
Their result addresses this question for lobsters but does not by
itself give a graceful supertree of the same size. The following
strengthening remains a natural question for \emph{general} trees:
\begin{equation}\label{eq:question-factor}
 s(T)\le(2-\varepsilon)m+K
 \quad\text{for absolute constants $\varepsilon>0$, $K$}.
\end{equation}
The inequality $\operatorname{grac}(T)\le s(T)$ shows that
\eqref{eq:question-factor} would answer van Bussel's question in that
generality. One may also ask whether the range-relaxed lobster
labelings of \cite{BK17} admit completion to graceful supertrees
without enlarging their label interval.

For each $\varepsilon>0$, Adamaszek, Allen, Grosu and Hladk\'y
\cite{AAGH20} proved that all sufficiently large $n$-vertex trees of
maximum degree at most $\eta n/\log n$, for a suitable $\eta>0$,
admit range-relaxed labels in
$\{1,\ldots,\lceil(1+\varepsilon)n\rceil\}$.
Erd\H{o}s, Hell and Winkler introduced \emph{gracesize} as the
maximum number of distinct differences induced by a bijection of
$V(T)$ to $\{1,\ldots,n\}$ \cite{EHW88}; see also
Rosa and \v{S}ir\'a\v{n} \cite{RS95}.
Letzter, Pokrovskiy and Williams \cite{LPW} established that every
sufficiently large $n$-vertex tree has gracesize at least
$(1-\varepsilon)n$. They distinguish this
notion from the all-tree range-relaxed question. In particular, a
uniform estimate $s(T)=m+o(m)$ would imply range-relaxed labelings of
span $m+o(m)$ for every tree, a conclusion not supplied by the
asymptotic gracesize result.

For perspective, trees with at most 35 vertices have been verified
graceful by Fang \cite{Fang}, and all trees of diameter at most five
are graceful \cite{HH01} (the latter includes the earlier results for
smaller diameter). On those classes $s(T)=m$. The purpose of the
bounds proved here is a uniform construction for arbitrary trees; in
particular, $s(T)=m$ for all $T$ is equivalent to the Graceful Tree
Conjecture.

\section*{Acknowledgments}
Several AI models were used in exploring constructions, checking finite
cases, and assisting with proof development. The author has verified the results and
takes responsibility for them.

\end{document}